\documentclass[12pt]{amsart}
\usepackage{mathrsfs}
\usepackage{amsfonts}
\usepackage{amssymb}
\usepackage[papersize={8.2in,11in},textwidth=16.4cm,textheight=22cm,centering]{geometry}
\usepackage{enumerate}

\usepackage{xcolor}
\definecolor{cite}{rgb}{0.30,0.60,1.00}
\definecolor{url}{rgb}{1.00,0.10,0.80}
\definecolor{link}{rgb}{0.00,0.00,1.00}

\usepackage[colorlinks,linkcolor=link,urlcolor=url,citecolor=cite,pagebackref,breaklinks]{hyperref}
\hypersetup{
pdfstartpage=1,
pdfstartview=FitH}

\usepackage{graphicx}

\usepackage{color}

\usepackage[OT2,T1]{fontenc}

\DeclareSymbolFont{cyrletters}{OT2}{wncyr}{m}{n}
\DeclareMathSymbol{\Sha}{\mathalpha}{cyrletters}{"58}

\DeclareFontFamily{U}{mathx}{\hyphenchar\font45}
\DeclareFontShape{U}{mathx}{m}{n}{
      <5> <6> <7> <8> <9> <10>
      <10.95> <12> <14.4> <17.28> <20.74> <24.88>
      mathx10
      }{}
\DeclareSymbolFont{mathx}{U}{mathx}{m}{n}
\DeclareMathAccent{\widecheck}{\mathalpha}{mathx}{"71}

 \usepackage{caption} 
\numberwithin{equation}{section}

\allowdisplaybreaks

\newtheorem*{theorem*}{Theorem}
\newtheorem{theorem}{Theorem}[section]
\newtheorem{lemma}{Lemma}[section]

\newtheorem{proposition}{Proposition}[section]

\newtheorem{corollary}{Corollary}[section]

\newtheorem*{claim*}{Claim}

\theoremstyle{remark}
\newtheorem{remark}{\bf Remark}
\newcommand{\ud}{\mathrm{d}}
\newcommand{\ue}{\mathrm{e}}

\newcommand{\kl}{\mathrm{Kl}}

\newcommand{\ind}{\mathrm{ind}}
\newcommand{\RRe}{\mathrm{Re}}
\newcommand{\IIm}{\mathrm{Im}}

\DeclareMathOperator{\Mod}{mod}

\renewcommand{\bmod}[1]{\,(\Mod{ #1})}

\newcommand{\bx}{\mathbf{x}}
\newcommand{\by}{\mathbf{y}}

\newcommand{\bC}{\mathbf{C}}
\newcommand{\bF}{\mathbf{F}}

\newcommand{\bR}{\mathbf{R}}
\newcommand{\bZ}{\mathbf{Z}}

\newcommand{\cC}{\mathcal{C}}
\newcommand{\cD}{\mathcal{D}}
\newcommand{\cE}{\mathcal{E}}
\newcommand{\cF}{\mathcal{F}}

\newcommand{\cH}{\mathcal{H}}

\newcommand{\cK}{\mathcal{K}}
\newcommand{\cL}{\mathcal{L}}
\newcommand{\cM}{\mathcal{M}}

\newcommand{\cR}{\mathcal{R}}

\newcommand{\fa}{\mathfrak{a}}
\newcommand{\fc}{\mathfrak{c}}

\newcommand{\fA}{\mathfrak{A}}

\newcommand{\fS}{\mathfrak{S}}

\begin{document}

\title[Multiplicative analogues of Kloosterman sums]{Moments and equidistributions of multiplicative analogues of Kloosterman sums}
\author{Ping Xi}

\address{School of Mathematics and Statistics, Xi'an Jiaotong University, Xi'an 710049, CHINA}
\email{ping.xi@xjtu.edu.cn}

\subjclass[2020]{11L40, 11M06, 11L05, 11T23, 60E05, 60E10}

\keywords{character sum, Kloosterman sum, Dirichlet $L$-function, moment, equidistribution, arcsine law}

\begin{abstract}
We consider a family of character sums as multiplicative analogues of Kloosterman sums. Using Gauss sums, Jacobi sums and Katz's bound for hypergeometric sums, we establish asymptotic formulae for any real (positive) moments of the above character sum as the character runs over all non-trivial multiplicative characters mod $p.$ Moreover, an arcsine law is also established as a consequence of the method of moments.
The evaluations of these moments also allow us to obtain asymptotic formulae for moments of such character sums weighted by special $L$-values (at $1/2$ and $1$). 
\end{abstract}

\dedicatory{{\it \small Dedicated to Professor Jie WU \\ on the occasion of his 61st birthday}}

\maketitle


\section{Introduction and the main results}
\subsection{Notation and conventions}
Throughout this paper, $p$ is always reserved for odd prime numbers. We do not distinguish $\bF_p$ with the complete residue system mod $p$, where the former one denotes the finite field of $p$ elements.

\begin{itemize}
\item For an integer $a$ coprime to $p,$ denote by $\overline{a}$ the multiplicative inverse of $a\bmod p$, i.e., $a\overline{a}\equiv1\bmod p.$
\item The symbol $*$ means to sum over primitive elements; e.g., primitive residue classes, or primitive characters as will be indicated.
\item Denote by $\phi$ the quadratic multiplicative character mod $p$, i.e., the Legendre symbol.
\item Put $\varepsilon_p=1$ if $p\equiv1\bmod4,$ and $=i$ if $p\equiv3\bmod4.$
\item We use $\varepsilon$ to denote a very small positive number, which might be different at each occurrence; we also write $X^\varepsilon \log X\ll X^\varepsilon.$ All implied constants in $O$ are absolute, with some exceptions that may depend on $\varepsilon$ only.
\end{itemize}

\subsection{Backgrounds}
Let $\chi$ be a non-trivial multiplicative character mod $p.$ The main concern here is the character sum
\begin{align}\label{eq:K}
K(\chi)=\frac{1}{\sqrt{p}}~\sideset{}{^*}\sum_{a\bmod p}\chi(a+\overline{a}).
\end{align}
This is clearly motivated by the (normalized) classical Kloosterman sum
\begin{align}\label{eq:Kloosterman}
\kl(\psi)=\frac{1}{\sqrt{p}}~\sideset{}{^*}\sum_{a\bmod p}\psi(a+\overline{a}),
\end{align}
where $\psi$ is chosen as a non-trivial additive character mod $p.$

The study of Kloosterman sums $\kl(\psi)$ has received considerable attention in past decades. As an immediate consequence of his proof on the Riemann Hypothesis for curves over finite fields, Weil \cite{We48} obtained the inequality
\begin{align*}
|\kl(\psi)|\leqslant2.
\end{align*}
Note that $\kl(\psi)$ is always real-valued, and its distribution in $[-2,2]$ becomes clear thanks to the work of Deligne and Katz. For instance, Katz \cite{Ka88} proved (with suitable generalizations) that the set
\begin{align*}
\{\kl(\psi):\psi\bmod p,\psi\text{~non-trivial}\}\subseteq[-2,2]
\end{align*}
becomes equidistributed with respect to the Sato--Tate measure $\frac{1}{2\pi}\sqrt{4-x^2}\ud x$ as $p\rightarrow+\infty$ over primes.

Different from Kloosterman sums, $K(\chi)$ is not real in general. It is obvious that $K(\chi)$ vanishes if $\chi$ is odd, i.e., $\chi(-1)=-1,$ and it follows from Weil \cite{We48} that
\begin{align*}
|K(\chi)|\leqslant2
\end{align*}
when $\chi$ is even and non-trivial. It is also natural to expect an equidistribution of $K(\chi)$  in the disk $\{z\in\bC:|z|\leqslant2\}$, or that of $|K(\chi)|$ in the interval $[-2,2]$,
as long as $\chi$ runs over a suitable family of even (non-trivial) multiplicative characters mod $p$. We will see this is in fact the case when all non-trivial even characters are taken into account, and the underlying idea is to consider the moment
\begin{align*}
M_k(p)=\sideset{}{^*}\sum_{\chi\bmod p}|K(\chi)|^{2k}.
\end{align*}
More generally, we put
\begin{align}\label{eq:moment}
M_\kappa(n,p)=\sideset{}{^*}\sum_{\chi\bmod p}\chi(n)|K(\chi)|^{2\kappa}
\end{align}
for any positive real number $\kappa$ and $n\in\bZ$.

\subsection{Equidistributions}
Since $K(\chi)$ is complex-valued in general, it should be easier to consider the distribution of the modulus $|K(\chi)|$ for the moment.
To this end, we henceforth write
\begin{align*}
\fA(p):=\{\chi\bmod p:~\chi(-1)=1,~\chi\text{~non-trivial}\}.
\end{align*}

\begin{theorem}\label{thm:|K|-equidistribution}
As $p\rightarrow+\infty$ over primes, the set $\{|K(\chi)|:\chi\in\fA(p)\}
$
becomes equidistributed in $[0,2]$ with respect to the measure
$\frac{2}{\pi}\frac{\ud x}{\sqrt{4-x^2}}.$
More precisely, for any interval $[a,b]\subseteq[0,2],$ we have
\begin{align*}
\lim_{p\rightarrow+\infty}\frac{|\{\chi\in\fA(p):|K(\chi)|\in[a,b]\}|}{(p-3)/2}=\frac{2}{\pi}(\arcsin(b/2)-\arcsin(a/2)).
\end{align*}
\end{theorem}

The measure
$\frac{2}{\pi}\frac{\ud x}{\sqrt{4-x^2}}$ reminds us the arcsine distribution, the probability density function of which is $\frac{1}{\pi\sqrt{x(1-x)}}.$
With a slight modification on the random variable $\chi\mapsto|K(\chi)|,$ we may conclude the following arcsine law from Theorem \ref{thm:|K|-equidistribution}.
\begin{corollary}\label{coro:|K|-arcsine}
As $p\rightarrow+\infty$ over primes, the random variables
\begin{align*}
\chi\mapsto\sqrt{|K(\chi)|/2}
\end{align*}
on $\fA(p)$ converge in distribution to the arcsine law. More precisely, for any given $u\in[0,1],$ we have
\begin{align*}
\lim_{p\rightarrow+\infty}\frac{|\{\chi\in\fA(p):\sqrt{|K(\chi)|/2}\leqslant u\}|}{(p-3)/2}=\frac{2}{\pi}\arcsin\sqrt{u}.
\end{align*}
\end{corollary}

In probability theory, the arcsine laws can appear in one-dimensional random walks and Brownian motion (the Wiener process). There are also some number theoretic models that have arcsine distributions (as special cases of Beta distributions). A very early example is due to Erd\H{o}s \cite{Er69}, asserting that
\begin{align*}
\lim_{x\rightarrow+\infty}\frac{1}{x}|\{n\leqslant x:s(n)\leqslant u\log\log n\}|=\frac{2}{\pi}\arcsin\sqrt{u},
\end{align*}
where $s(n)$ denotes the number of $j\leqslant \omega(n)$ such that $j$-th prime factor $p_j$ of $n$ (in increasing order and counting without multiplicities) satisfies $\log\log p_j<j.$
Another instance is the so-called DDT Theorem proven by Deshouillers, Dress and Tenenbaum \cite{DDT79} concerning the random variable $D_n$ which takes the value $\log d$ with the uniform probability $1/d(n)$, where $d(n)$ counts the number of (positive) divisors of $n$. Note that the expectation of $D_n$ is $\frac{1}{2}\log n$, and it is proven in \cite{DDT79}, as $x\rightarrow+\infty,$ that
\begin{align*}
\frac{1}{x}\sum_{n\leqslant x}\mathrm{Prob}(D_n\leqslant u\log n)=\frac{2}{\pi}\arcsin\sqrt{u}+O((\log x)^{-\frac{1}{2}})
\end{align*}
uniformly in $u\in[0,1].$ See \cite{Te80,BM07,CW14,BT16,DT18} for instance with many other extensions and further developments. Corollary \ref{coro:|K|-arcsine} above provides another interesting example that also has the arcsine law.

Since $K(\chi)$ may take complex values, it is expected that Theorem \ref{thm:|K|-equidistribution} can be generalized to a two-dimensional equidistribution of $K(\chi)$ in the disk $\{z\in\bC:|z|\leqslant2\}$. 

\begin{theorem}\label{thm:K-equidistribution}
As $p\rightarrow+\infty$ over primes, the set $\{K(\chi):\chi\in\fA(p)\}$
becomes equidistributed in $\cD:=\{x+iy: x^2+y^2\leqslant4,x,y\in\bR\}$ with respect to the measure
\begin{align*}
\ud\mu&:=\frac{1}{\pi^2}\frac{\ud x\ud y}{\sqrt{x^2+y^2}\sqrt{4-(x^2+y^2)}}.
\end{align*}
More precisely, for any domain $\cR\subseteq\cD,$ we have
\begin{align*}
\lim_{p\rightarrow+\infty}\frac{|\{\chi\in\fA(p):K(\chi)\in \cR\}|}{(p-3)/2}=\iint_{\cR} \ud\mu.
\end{align*}
\end{theorem}

Making the change of variables $(x,y)\rightarrow (r\cos\theta,r\sin\theta)$, we then find $\ud\mu$ is changed to
\begin{align*}
\frac{1}{\pi^2}\frac{\ud r\ud\theta}{\sqrt{4-r^2}}=\frac{2}{\pi}\frac{\ud r}{\sqrt{4-r^2}}\cdot\frac{\ud\theta}{2\pi},
\end{align*}
which is exactly the product of the Lebesgue measure $\frac{1}{2\pi}\ud\theta$ and the probability measure in Theorem \ref{thm:|K|-equidistribution}. Since $r$ stands for the modulus in $[0,2]$ and $\theta$ characterizes the arguments in $[-\pi,\pi]$, we may conclude the following equidistribution in the polar coordinates.
\begin{theorem}\label{thm:K-equidistribution2}
For any intervals $I\subseteq[-\pi,\pi]$ and $J\subseteq[0,2],$ we have
\begin{align*}
\lim_{p\rightarrow+\infty}\frac{|\{\chi\in\fA(p):\arg K(\chi)\in I,|K(\chi)|\in J\}|}{(p-3)/2}=\frac{|I|}{\pi^2}\int_J\frac{\ud r}{\sqrt{4-r^2}}.
\end{align*}
\end{theorem}

Theorem \ref{thm:|K|-equidistribution} is now an immediately consequence of Theorem \ref{thm:K-equidistribution2} by taking
$I=[-\pi,\pi].$ On the other hand, taking $J=[0,2]$ in Theorem \ref{thm:K-equidistribution2}, we can also obtain a uniform distribution of $\arg K(\chi)$ in $[-\pi,\pi].$ Therefore, as one might have imagined, $\{K(\chi):\chi\in\fA(p)\}$ generate two random variables $\chi\mapsto\arg K(\chi)$ and $\chi\mapsto|K(\chi)|$ which become independent as $p\rightarrow+\infty$.

We will give a quick proof of Theorem \ref{thm:|K|-equidistribution} in Section \ref{sec:|K|-equidistribution-proof} using the method of moments (see Appendix \ref{sec:methodofmoments}), although we may conclude it from Theorem \ref{thm:K-equidistribution2}. The proof therein will be helpful to understand the proof of the latter one in 
Section \ref{sec:K-equidistribution-proof}.

Although $K(\chi)$ is complex-valued in general, we will see from Lemma \ref{lm:K-transform} that the normalized sum $\varepsilon_p\overline{\chi}(2)K(\chi)$ is always real-valued and falls into $[-2,2]$ for any non-trivial $\chi\bmod p$ with $\chi(-1)=1.$ Based on a similar argument in proving Theorems \ref{thm:|K|-equidistribution} and \ref{thm:K-equidistribution}, we may conclude the following equidistribution.

\begin{theorem}\label{thm:K-twistedequidistribution}
As $p\rightarrow+\infty$ over primes, the set
\begin{align*}
\{\varepsilon_p\overline{\chi}(2)K(\chi):\chi\in\fA(p)\}\subseteq[-2,2]
\end{align*}
becomes equidistributed with respect to the measure $\frac{1}{\pi}\frac{\ud x}{\sqrt{4-x^2}}.$ More precisely, for any interval $[a,b]\subseteq[-2,2],$ we have
\begin{align*}
\lim_{p\rightarrow+\infty}\frac{|\{\chi\in\fA(p):\varepsilon_p\overline{\chi}(2)K(\chi)\in[a,b]\}|}{(p-3)/2}=\frac{1}{\pi}(\arcsin(b/2)-\arcsin(a/2)).
\end{align*}
\end{theorem}

\subsection{Moments}
As mentioned above, the underlying idea to prove equidistributions is evaluating the moment $M_\kappa(n,p)$ asymptotically as $p\rightarrow+\infty$.

\begin{theorem}\label{thm:K-moment}
Let $\kappa$ be a fixed positive real number. For all large primes $p,$ we have
\begin{align*}
M_\kappa(n,p)=\frac{\Gamma(2\kappa+1)}{2\Gamma(\kappa+1)^2}\cdot \Delta_n\cdot p+O(\kappa4^{\kappa}p^{\frac{1}{2}}\log p),
\end{align*}
where $\Delta_n=1$ or $0$ according to $n\equiv\pm1\bmod p$ or not, and the implied constant is absolute.
\end{theorem}

\begin{corollary}\label{coro:K-moment}
Let $k$ be a fixed positive integer. For all large primes $p,$ we have
\begin{align*}
M_k(n,p)=\frac{1}{2} \binom{2k}{k}\cdot \Delta_n\cdot p+O(k4^kp^{\frac{1}{2}}\log p),
\end{align*}
where $\Delta_n$ is defined as above and the implied constant is absolute.
\end{corollary}

Thanks to the method of moments, one can imagine that Theorem \ref{thm:|K|-equidistribution} should be a direct consequence of Theorem \ref{thm:K-moment}. On the other hand, the proof of Theorem  \ref{thm:K-equidistribution} would require to study the mixed moment
\begin{align}\label{eq:mixedmoment}
M_{k,l}(n,p)=\sideset{}{^*}\sum_{\chi\bmod p}\chi(n)\overline{K(\chi)}^k K(\chi)^l
\end{align}
for all integers $k,l\geqslant0$ and $n\in\bZ$.
Clearly, we have $M_{k,k}(n,p)=M_k(n,p).$ In general, we have following asymptotic evaluations.

\begin{theorem}\label{thm:K-mixedmoment}
Let $k,l$ be two non-negative integers. For all large primes $p,$ we have
\begin{align*}
M_{k,l}(n,p)
&=\varepsilon_p^{k-l}\binom{k+l}{(k+l)/2}\cdot \frac{\Delta_{n;k,l}}{2}\cdot p+O((k+l)2^{k+l}p^{\frac{1}{2}}\log p),
\end{align*}
where $\Delta_{n;k,l}=1$ or $0$ according to $2^k\equiv\pm2^ln\bmod p$ or not, and the implied constant is absolute.
\end{theorem}

As another consequence of Theorem \ref{thm:K-moment}, we can conclude asymptotic formulae for moments of $K(\chi)$ twisted by Dirichlet $L$-functions. To this end, let $\beta$ be the multiplicative function defined by $\beta(p^\nu)=4^{-\nu}\binom{2\nu}{\nu}$ for all primes $p$ and $\nu\in\bZ^+.$ Define 
\[\fS:=\sum_{n\geqslant1}\frac{\beta(n)^2}{n^2}.\]

\begin{corollary}\label{coro:K-moment-L1}
Let $\kappa$ be a fixed positive real number. For all large primes $p,$ we have
\begin{align*}
\sideset{}{^*}\sum_{\chi\bmod p}|K(\chi)|^{2\kappa}|L(1,\chi)|=\frac{\Gamma(2\kappa+1)}{2\Gamma(\kappa+1)^2}\cdot \fS\cdot p+O(4^{\kappa}p^{\frac{1}{2}}\log p),
\end{align*}
where the implied constant is absolute.
\end{corollary}

The special case $n=\pm1$ in Corollary \ref{coro:K-moment} confirms a conjecture of Zhang \cite{BRZ22}.
Corollary \ref{coro:K-moment-L1} solves Conjecture 1.11 in Bag, Rojas-Le\'on and Zhang \cite{BRZ22}, and thus solves the conjecture of Zhang \cite{Zh02} in view of Corollary \ref{coro:K-moment}.

With more efforts, we can also characterize the moment of $K(\chi)$ weighted by $|L(1/2,\chi)|^2$.
\begin{theorem}\label{thm:K-moment-L1/2}
Let $\kappa$ be a fixed positive real number. For all large primes $p,$ we have
\begin{align*}
\sideset{}{^*}\sum_{\chi\bmod p}|K(\chi)|^{2\kappa}|L(1/2,\chi)|^2=\frac{\Gamma(2\kappa+1)}{2\Gamma(\kappa+1)^2}\cdot p\cdot (\log\frac{p}{8\pi}-\frac{\pi}{2}+\gamma)+O(5^\kappa p^{1-\delta})
\end{align*}
for any $\delta<1/8,$
where $\gamma$ is the Euler constant and the implied constant depends only on $\delta$.
\end{theorem}

The proof of Theorem \ref{thm:K-moment-L1/2} is also based on that of Theorem \ref{thm:K-moment}, expressing this weighted moment in terms of $M_\kappa(n,p)$ with suitable $n$. The main difficulty would come from bounding the error term, for which we will appeal to a deep estimate for smooth bilinear forms due to Fouvry, Kowalski and Michel \cite{FKM14}, and thus more algebro-geometric tools come into the picture.

The character sum $K(\chi)$ is ultimately related to the generalized Gauss sum
\[\tau_2(n,\chi)=\sum_{v\bmod p}\chi(v)\ue\Big(\frac{nv^2}{p}\Big).\]
Zhang \cite{Zh02} studied the $2k$-th moment of $|\tau_2(n,\chi)|$ for $k=2,3,4$ as $\chi$ runs over all non-trivial characters mod $p$. Following the same arguments as in proving Theorem \ref{thm:K-moment}, we may also obtain all higher moments of $|\tau_2(n,\chi)|$, even with twists by Dirichlet $L$-functions, and we will come back to this issue in the last section. 

It is worthy to point out that the work in this paper is partly motivated by Zhang \cite{Zh02}, Bag and Barman \cite{BB21} and Bag, Rojas-Le\'on and Zhang \cite{BRZ22}.
Bag and Barman \cite{BB21} obtained asymptotic formulae for $2k$-th moments of $\tau_2(n,\chi)$ and $K(\chi)$  for some small integral $k$.

Our arguments are quite different from those in  \cite{BB21} and \cite{BRZ22}, and the tools here will include Gauss sums, Jacobi sums and also hypergeometric sums as we will see.

Motived by hyper-Kloosterman sums (see \eqref{eq:hyperKloosterman} below), one may also consider the high-dimensional generalization of $K(\chi):$
\begin{align}\label{eq:K-hyper}
K_k(a,\chi)=p^{\frac{1-k}{2}}~\mathop{\sum\cdots\sum}_{\substack{a_1,a_2,\cdots,a_k\bmod p\\ a_1a_2\cdots a_k\equiv a\bmod p}}\chi(a_1+a_2+\cdots+a_k).
\end{align}
Clearly, $K(\chi)=K_2(1,\chi).$ The arguments in this paper would also work for studying the distribution and moments of $K_k(a,\chi)$, and we will address this issue in Section \ref{sec:remarks}.

\subsection*{Acknowledgements}
I would like to thank Professor Wenpeng Zhang for introducing me to his conjectures. Sincere thanks are also due to the referee for his/her kind suggestion to conclude equidistributions.
The work is supported in part by NSFC (No. 12025106, No. 11971370).

\smallskip

\section{Preliminaries}
Many of the statements in this section are standard, and the readers can refer to many classical monographs (e.g., \cite{IK04}).

\subsection{Gauss sums} 
Given a multiplicative character $\chi\bmod p$ and $n\in\bZ,$ the Gauss sum is defined by
\begin{align*}
\tau(n,\chi)=\sum_{v\bmod p}\chi(v)\ue\Big(\frac{nv}{p}\Big).
\end{align*}
In what follows, we write $\tau(\chi)=\tau(1,\chi).$ It is well-known that
$|\tau(\chi)|=\sqrt{p}$ when $\chi$ is non-trivial. Recall that $\phi$ denotes the quadratic character mod $p$. The evaluation of $\tau(n,\phi)$ can be made explicitly as follows.

\begin{lemma}\label{lm:Gauss-evaluation}
For each integer $n,$ we have
\begin{align*}
\tau(n,\phi)=\varepsilon_p\phi(n)\sqrt{p},
\end{align*}
where $\varepsilon_p=1$ if $p\equiv1\bmod4,$ and $=i$ if $p\equiv3\bmod4.$ In particular, $\tau(\phi)=\varepsilon_p\sqrt{p}.$
\end{lemma}

\subsection{Jacobi sums} 
We also use the Jacobi sum
\begin{align}\label{eq:Jacobisum}
J(\chi_1,\chi_2)=\sum_{v\bmod p}\chi_1(v)\chi_2(1-v)
\end{align}
for multiplicative characters $\chi_1,\chi_2$ mod $p.$
In most cases, Jacobi sums can be expressed in terms of Gauss sums as follows.
\begin{lemma}\label{lm:Jacobi-Gauss}
Suppose $\chi_1,\chi_2$ and $\chi_1\chi_2$ are all non-trivial multiplicative characters.
Then
\begin{align*}
J(\chi_1,\chi_2)=\frac{\tau(\chi_1)\tau(\chi_2)}{\tau(\chi_1\chi_2)}.
\end{align*}
\end{lemma}

\subsection{Hypergeometric sums} 
We now consider hypergeometric sums introduced by Katz (see \cite[Chapter 8]{Ka90}), and the formulation is given in the settings of finite fields. Let $m,n$ be two non-negative integers. Suppose $\boldsymbol{\chi}=(\chi_i)_{1\leqslant i\leqslant m}$ and $\boldsymbol{\eta}=(\eta_j)_{1\leqslant j\leqslant n}$ are two tuples of multiplicative characters of $\bF_p^\times$ (which can be identified as Dirichlet characters modulo $p$ in a natural way). Katz introduced the following hypergeometric sum
\begin{align*}
H(t;\boldsymbol\chi,\boldsymbol\eta,p)
:=\frac{(-1)^{m+n-1}}{p^{(m+n-1)/2}}\mathop{\sum\sum}_{\substack{\bx\in(\bF_p^\times)^m,\by\in(\bF_p^\times)^n\\ N(\bx)=tN(\by)}}\boldsymbol\chi(\bx)
\overline{\boldsymbol\eta(\by)}\psi(T(\bx)-T(\by))
\end{align*}
for $t\in \bF_p^\times,$ where, for $\bx=(x_1,x_2,\cdots,x_m)\in(\bF_p^\times)^m$, 
\begin{align*}
\boldsymbol\chi(\bx)=\prod_{1\leqslant i\leqslant m}\chi_i(x_i),
\end{align*}
\begin{align}\label{eq:T(x)N(x)}
T(\bx)=x_1+x_2+\cdots+x_m,\ \ \ N(\bx)=x_1x_2\cdots x_m,
\end{align}
and the notation with $\by$ can be defined in the same way. We say $\boldsymbol{\chi}$ and $\boldsymbol{\eta}$ are {\it disjoint} if $\chi_i\neq\eta_j$ for all $1\leqslant i\leqslant m$ and $1\leqslant j\leqslant n.$

As one may see, $H(t;\boldsymbol\chi,\boldsymbol\eta,p)$ extends the classical (hyper-) Kloosterman sums. In particular, if $n=0$ and $\chi_i=\mathbf{1}$ for each $1\leqslant i\leqslant m,$ the hypergeometric sum becomes the hyper-Kloosterman sum
\begin{align}\label{eq:hyperKloosterman}
\kl_m(t,p)=\frac{(-1)^{m-1}}{q^{(m-1)/2}}\sum_{\substack{\bx\in(\bF_p^\times)^m\\ N(\bx)=t}}\psi(T(\bx)),
\end{align}
which was already proved by Deligne \cite{De80} that $t\mapsto\kl_m(t,p)$ is the Frobenius trace function of a certain $\ell$-adic sheaf ($\ell\neq p$), which is geometrically irreducible of rank $m$, lisse on $\mathbf{G}_{m,\bF_p}$ and pointwise pure of weight 0.
In general, Katz \cite[Theorem 8.4.2]{Ka90} proved the following assertion, which is quite crucial in proving equidistributions as required in this paper.

\begin{lemma}\label{lm:hypergeometricsum}
With the above notation,
if $\boldsymbol\chi$ and $\boldsymbol\eta$ are disjoint, then for any $\ell\neq p,$ there exists a geometrically irreducible $\ell$-adic middle-extension sheaf $\cH(\boldsymbol\chi,\boldsymbol\eta)$ on $\mathbf{A}_{k}^1$ with trace function given by $t\mapsto H(t;\boldsymbol\chi,\boldsymbol\eta,p),$ such that it is
\begin{itemize}
\item pointwise pure of weight $0$ and of rank $\max\{m,n\};$
\item lisse on $\mathbf{G}_{m,\bF_p}$, if $m\neq n;$
\item lisse on $\mathbf{G}_{m,\bF_p}-\{1\}$ and of rank $m,$ if $m=n.$
\end{itemize}
\end{lemma}

Trivially we have
\begin{align*}
|H(t;\boldsymbol\chi,\boldsymbol\eta,p)|\leqslant p^{(m+n-1)/2}.
\end{align*}
As a direct consequence, in the case that $\boldsymbol\chi$ and $\boldsymbol\eta$ are disjoint, we can bound the hypergeometric sum as
\begin{align*}
|H(t;\boldsymbol\chi,\boldsymbol\eta,p)|\leqslant \max\{m,n\},
\end{align*}
illustrating squareroot cancellations within the summations.

\subsection{Transformations of $K(\chi)$}
The following lemma expresses $K(\chi)$ as nonlinear character sums or Gauss sums.
\begin{lemma}\label{lm:K-transform}
Let $\chi$ be a non-trivial multiplicative character mod $p$ with $\chi(-1)=1.$ Then
\begin{align*}
K(\chi)=\frac{\chi(2)}{\sqrt{p}}\sum_{b\bmod p}\chi(b)\phi(b^2-1).
\end{align*}
Moreover, suppose $\chi=\chi_1^2$ for some character $\chi_1,$ then we have
\begin{align}\label{eq:K-Gausssum}
K(\chi)=\frac{2\chi(2)\overline{\tau(\phi)}}{p^{3/2}}\cdot \RRe(\tau(\chi_1)\overline{\tau(\chi_1\phi)})
=\frac{2\chi(2)\overline{\varepsilon_p}}{p}\cdot \RRe(\tau(\chi_1)\overline{\tau(\chi_1\phi)}).
\end{align}
\end{lemma}

\proof
We first write
\begin{align*}
K(\chi)=\frac{1}{\sqrt{p}}\sum_{b\bmod p}\chi(b)\sum_{\substack{a\bmod p\\ a+\overline{a}\equiv b\bmod p}}1=\frac{1}{\sqrt{p}}\sum_{b\bmod p}\chi(b)\sum_{\substack{a\bmod p\\ a^2-ba+1\equiv 0\bmod p}}1.
\end{align*}
The first part follows by noting that the number of solutions to $a^2-ba+1\equiv 0\bmod p$ in $a\bmod p$ is exactly $1+\phi(b^2-4).$ 

Regarding the second part, we write
\begin{align*}
K(\chi)
&=\frac{\chi(2)}{\sqrt{p}}\sum_{b\bmod p}\chi_1(b^2)\phi(b^2-1)\\
&=\frac{\chi(2)}{\sqrt{p}}\sum_{b\bmod p}(1+\phi(b))\chi_1(b)\phi(b-1).
\end{align*}
Therefore,
\begin{align*}
K(\chi)=\frac{\chi(2)\phi(-1)}{\sqrt{p}}(J(\chi_1,\phi)+J(\chi_1\phi,\phi)).
\end{align*}
Since $\chi$ is non-trivial, we must have $\chi_1\neq\phi$, and it then follows from Lemma \ref{lm:Jacobi-Gauss} that
\begin{align*}
K(\chi)
&=\frac{\chi(2)\phi(-1)}{\sqrt{p}}\Big(\frac{\tau(\chi_1)\tau(\phi)}{\tau(\chi_1\phi)}+\frac{\tau(\chi_1\phi)\tau(\phi)}{\tau(\chi_1)}\Big)\\
&=\frac{\chi(2)\phi(-1)\tau(\phi)}{p^{3/2}}(\tau(\chi_1)\overline{\tau(\chi_1\phi)}+\tau(\chi_1\phi)\overline{\tau(\chi_1)}),
\end{align*}
which completes the proof readily in view of $\phi(-1)\tau(\phi)=\overline{\tau(\phi)}$.
\endproof

\begin{lemma}\label{lm:twoGauss-moments}
Let $k$ be a fixed positive integer and $(n,p)=1$. Then we have
\begin{align*}
\sum_{\chi\bmod p}\overline{\chi}(n)(\tau(\chi)\overline{\tau(\chi\phi)})^k
=-(p-1)p^{k-\frac{1}{2}}H(n;\mathbf{1},\boldsymbol\phi,p),
\end{align*}
where $\boldsymbol\phi=(\phi,\phi,\cdots,\phi)$ is the $k$-tuple of characters given by $\phi\bmod p.$
\end{lemma}

\proof
Denote by $T$ the sum in question. By definition we may write
\begin{align*}
T&=\mathop{\sum\sum}_{\bx,\by\in[1,p[^k}\phi(N(\by))\ue\Big(\frac{T(\bx)-T(\by)}{p}\Big)\sum_{\chi\bmod p}\chi(N(\bx))\overline{\chi}(nN(\by)),
\end{align*}
where, for $\bx\in[1,p[^k,$ $T(\bx)$ and $N(\bx)$ are defined the same with \eqref{eq:T(x)N(x)}.
The lemma follows immediately from the orthogonality of multiplicative characters.
\endproof

\begin{lemma}\label{lm:Jacobi-moments}
Let $k$ be a fixed positive integer and $(n,p)=1$. Then we have
\begin{align*}
\sum_{\chi\bmod p}\overline{\chi}(n)J(\chi,\phi)^k=-(p-1)p^{-\frac{1}{2}}\tau(\phi)^kH(n;\mathbf{1},\boldsymbol\phi,p).
\end{align*}
Moreover, we have
\begin{align*}
\Big|\sum_{\chi\bmod p}\overline{\chi}(n)J(\chi,\phi)^k\Big|\leqslant kp^{\frac{k+1}{2}}.
\end{align*}
\end{lemma}

\proof Denote by $T$ the sum in question. Invoking Lemma \ref{lm:Jacobi-Gauss}, we may write
\begin{align*}
T&=\sum_{\substack{\chi\bmod p\\ \chi\neq \chi_0,\phi}}\overline{\chi}(n)J(\chi,\phi)^k+\phi(n)(-\phi(-1))^k+(-1)^k\\
&=\sum_{\substack{\chi\bmod p\\ \chi\neq \chi_0,\phi}}\overline{\chi}(n)\Big(\frac{\tau(\chi)\tau(\phi)}{\tau(\chi\phi)}\Big)^k+\phi(n)(-\phi(-1))^k+(-1)^k\\
&=p^{-k}\tau(\phi)^k\sum_{\substack{\chi\bmod p\\ \chi\neq \chi_0,\phi}}\overline{\chi}(n)(\tau(\chi)\overline{\tau(\chi\phi)})^k+\phi(n)(-\phi(-1))^k+(-1)^k.\end{align*}
The first part of the lemma then follows by noting that $\tau(\phi)^2=\phi(-1)p,$ $\tau(\chi_0)=-1$
and Lemma \ref{lm:twoGauss-moments}. And the second part is an immediate consequence of Lemma \ref{lm:hypergeometricsum}.
\endproof

\smallskip

\section{Group of characters}\label{sec:characters}

Denote by $\cC_p$ be the set of all multiplicative characters mod $p$. This forms a cyclic group of order $p-1$, which is isomorphic to $(\bZ/p\bZ)^\times.$ Fix $\eta$ as a generator in $\cC_p.$ It is clearly that $\eta(a)=1$ if and only if $a\equiv1\bmod p.$
Denote by $\cC_p^+$ the set of all even characters. This is an index two subgroup of $\cC_p$, and we may write
$$\cC_p^+=\{\eta^{2i}:1\leqslant i\leqslant (p-1)/2\}.$$

\begin{lemma}\label{lm:Gaussmoment-j}
Let $k$ be a fixed positive integer and $(a,p)=1$. Then we have
\begin{align*}
\sum_{1\leqslant i\leqslant \frac{p-1}{2}}\eta^i(a)(\tau(\eta^i)\overline{\tau(\eta^i\phi)})^k&\ll kp^{k+\frac{1}{2}}\log p,
\end{align*}
where the implied constant is absolute.
\end{lemma}

\proof
Denote by $\Sigma$ the sum in question. By definition, we may write
\begin{align*}
\Sigma
&=\mathop{\sum\sum}_{\bx,\by\in[1,p[^k}\phi(N(\by))\ue\Big(\frac{T(\bx)-T(\by)}{p}\Big)\sum_{1\leqslant i\leqslant \frac{p-1}{2}}\eta(aN(\bx)/N(\by))^i\\
&=\mathop{\sum\sum}_{\bx,\by\in[1,p[^k}\phi(N(\by))\ue\Big(\frac{\bx\cdot \by-T(\by)}{p}\Big)\sum_{1\leqslant i\leqslant \frac{p-1}{2}}\eta(aN(\bx))^i\\
&=\frac{p-1}{2}\mathop{\sum\sum}_{\substack{\bx,\by\in[1,p[^k\\\eta(aN(\bx))=1}}\phi(N(\by))\ue\Big(\frac{\bx\cdot\by-T(\by)}{p}\Big)\\
&\ \ \ \ \ +\mathop{\sum\sum}_{\substack{\bx,\by\in[1,p[^k\\\eta(aN(\bx))\neq1}}\phi(N(\by))\ue\Big(\frac{\bx\cdot\by-T(\by)}{p}\Big)\frac{\eta(N(\bx))}{1-\eta(N(\bx))}(1-\phi(N(\bx))),
\end{align*}
where we also introduced suitable changes of variables.
By Lemma \ref{lm:Gauss-evaluation}, we have
\begin{align*}
\mathop{\sum}_{\by\in[1,p[^k}\phi(N(\by))\ue\Big(\frac{\bx\cdot\by-T(\by)}{p}\Big)
&=\prod_{1\leqslant j\leqslant k}\sum_{y_j\bmod p}\phi(y_j)\ue\Big(\frac{(x_j-1)y_j}{p}\Big)\\
&=(\varepsilon_p\sqrt{p})^k\prod_{1\leqslant j\leqslant k}\phi(x_j-1).
\end{align*}
It then follows that
\begin{align*}
\Sigma
&=(\varepsilon_p\sqrt{p})^k\Big(\frac{p-1}{2}\Sigma_1+\Sigma_2\Big),
\end{align*}
where
\begin{align*}
\Sigma_1&:=\mathop{\sum\cdots\sum}_{\substack{x_1\cdots,x_k\bmod p\\\eta(ax_1\cdots x_k)=1}}\phi((x_1-1)\cdots(x_k-1)),\\
\Sigma_2&:=\mathop{\sum\cdots\sum}_{\substack{x_1\cdots,x_k\bmod p\\\eta(ax_1\cdots x_k)\neq1}}\phi((x_1-1)\cdots(x_k-1))\frac{1-\phi(x_1\cdots x_k)}{1-\overline{\eta}(x_1\cdots x_k)}.
\end{align*}

First, we have
\begin{align*}
\Sigma_1
&=\frac{1}{p-1}\sum_{1\leqslant j\leqslant p-1}\mathop{\sum\cdots\sum}_{x_1\cdots,x_k\bmod p}\phi((x_1-1)\cdots(x_k-1))\eta(ax_1\cdots x_k)^j\\
&=\frac{\phi(-1)^k}{p-1}\sum_{1\leqslant j\leqslant p-1}\eta^j(a)J(\eta^j,\phi)^k.
\end{align*}
From Lemma \ref{lm:Jacobi-moments}, it follows that
\begin{align*}
\Sigma_1&\ll kp^{\frac{k-1}{2}}.
\end{align*}

We now turn to $\Sigma_2.$ Denote by $\mu_{p-1}$ the group of $(p-1)$-th roots of unity. Then
\begin{align*}
\Sigma_2
&=\sum_{\substack{\nu\in\mu_{p-1}\\\nu\neq1}}\frac{1}{1-\overline{\nu}}\mathop{\sum\cdots\sum}_{\substack{x_1\cdots,x_k\bmod p\\\eta(ax_1\cdots x_k)=\nu}}\phi((x_1-1)\cdots(x_k-1))(1-\phi(x_1\cdots x_k))\\
&=\frac{1}{p-1}\sum_{1\leqslant j\leqslant p-1}\sum_{\substack{\nu\in\mu_{p-1}\\\nu\neq1}}\frac{\nu^{-j}}{1-\overline{\nu}}\\
&\ \ \ \ \times\mathop{\sum\cdots\sum}_{x_1\cdots,x_k\bmod p}\phi((x_1-1)\cdots(x_k-1))(1-\phi(x_1\cdots x_k))\eta(ax_1\cdots x_k)^j\\
&=\frac{\phi(-1)^k}{p-1}\sum_{\substack{\nu\in\mu_{p-1}\\\nu\neq1}}\frac{1}{1-\overline{\nu}}\sum_{1\leqslant j\leqslant p-1}\nu^{-j}\eta^j(a)(J(\eta^j,\phi)^k-J(\eta^j\phi,\phi)^k).
\end{align*}

Given a multiplicative character $\xi\bmod p$, $\nu\in\mu_{p-1}$ and $(a,p)=1$, we consider
\begin{align*}
R_k(\nu,a,\xi):=\sum_{1\leqslant j\leqslant p-1}\nu^{-j}\eta^j(a)J(\eta^j\xi,\phi)^k,
\end{align*}
so that
\begin{align}\label{eq:Sigma2-initial}
\Sigma_2
&=\frac{\phi(-1)^k}{p-1}\sum_{\substack{\nu\in\mu_{p-1}\\\nu\neq1}}\frac{R_k(\nu,a,\chi_0)-R_k(\nu,a,\phi)}{1-\overline{\nu}}.
\end{align}

Without loss of generality, we assume $\nu=\ue(\frac{h}{p-1})$ for some $h\in \bZ/(p-1)\bZ$, so that
\begin{align*}
R_k(\nu,a,\xi)&=\sum_{1\leqslant j\leqslant p-1}\ue\Big(\frac{-jh}{p-1}\Big)\eta^j(a)J(\eta^j\xi,\phi)^k.
\end{align*}
Suppose $\eta$ is given by
\begin{align*}
n\mapsto\ue\Big(\frac{\ind_g(n)}{p-1}\Big),\ \ n\in(\bZ/p\bZ)^\times,
\end{align*}
where $g$ is a fixed primitive root mod $p$, and $\ind_g(n)$ is defined to be the index of $n$ $\bmod p$ with base $g$; i.e., the smallest non-negative  integer $\delta$ such that $g^\delta\equiv n\bmod p.$
Note that $\eta^j$ runs over all multiplicative characters mod $p$ as long as $j$ runs over all integers in $[1,p-1].$ Let $\chi$ be a multiplicative character mod $p$ given by
\[\chi(n)=\ue\Big(\frac{r\ind_g(n)}{p-1}\Big),\ \ n\in(\bZ/p\bZ)^\times\]
for some $r\in\bZ.$
If $\chi=\eta^j$ for some $1\leqslant j\leqslant p-1,$ then $r\equiv j\bmod{p-1}.$
In this way, we may write
\begin{align*}
\ue\Big(\frac{-jh}{p-1}\Big)
&=\ue\Big(\frac{-rh}{p-1}\Big)=\ue\Big(\frac{-r\ind_g(g^h)}{p-1}\Big)=\overline{\chi}(g^h).
\end{align*}
Therefore, we conclude from Lemma \ref{lm:Jacobi-moments} that
\begin{align*}
R_k(\nu,a,\xi)&=\sum_{\chi\bmod p}\overline{\chi}(g^h)\chi(a)J(\chi\xi,\phi)^k\\
&=\xi(g^h)\overline{\xi}(a)\sum_{\chi\bmod p}\overline{\chi}(g^h)\chi(a)J(\chi,\phi)^k.
\end{align*}
Substituting this into \eqref{eq:Sigma2-initial}, we derive from Lemma \ref{lm:hypergeometricsum} that
\begin{align*}
\Sigma_2
&\ll kp^{\frac{k-1}{2}}\sum_{1\leqslant |h|\leqslant (p-1)/2}\frac{1}{|1-\ue(\frac{h}{p-1})|}\\
&\ll kp^{\frac{k-1}{2}}\sum_{1\leqslant |h|\leqslant (p-1)/2}\frac{p}{|h|}\\
&\ll kp^{\frac{k+1}{2}}\log p.
\end{align*}

Combining the above estimates for $\Sigma_1$ and $\Sigma_2,$ we obtain
\begin{align*}
\Sigma
&\ll kp^{k+\frac{1}{2}}\log p,
\end{align*}
where the implied constant is absolute.
\endproof

\smallskip

\section{Proof of Theorem \ref{thm:K-moment}}\label{sec:proof-Theorem-K-moment}
As in Section \ref{sec:characters}, we fix $\eta$ as a generator of $\cC_p.$ It now follows from \eqref{eq:K-Gausssum} that
\begin{align*}
M_\kappa(n,p)
&=p^{-2\kappa}\sum_{1\leqslant i\leqslant \frac{p-3}{2}}\eta^{i}(n^2)\Big|\tau(\eta^i)\overline{\tau(\eta^i\phi)}+\overline{\tau(\eta^i)}\tau(\eta^i\phi)\Big|^{2\kappa}.
\end{align*}
For each $1\leqslant i\leqslant \frac{p-3}{2}$, we have
\begin{align*}
\Big|\tau(\eta^i)\overline{\tau(\eta^i\phi)}+\overline{\tau(\eta^i)}\tau(\eta^i\phi)\Big|^2=2p^2+2\Re(\tau(\eta^i)^2\overline{\tau(\eta^i\phi)}^2).
\end{align*}
Expanding the $\kappa$-th power, we may write
\begin{align}
M_\kappa(n,p)
&=(2/p^2)^\kappa\sum_{1\leqslant i\leqslant \frac{p-1}{2}}\eta^{i}(n^2)(p^2+\Re(\tau(\eta^i)^2\overline{\tau(\eta^i\phi)}^2))^\kappa+O(p^{-\kappa})\nonumber\\
&=p^{-2\kappa}\Gamma(\kappa+1)\sum_{j\geqslant 0}\frac{(2p^2)^{\kappa-j}}{j!\Gamma(\kappa-j+1)}N_j(\eta,n,p)+O(p^{-\kappa}),
\label{eq:Mk(p)-decomposition}
\end{align}
where
\begin{align*}
N_j(\eta,n,p)
&=\sum_{1\leqslant i\leqslant \frac{p-1}{2}}\eta^{i}(n^2)(\tau(\eta^i)^2\overline{\tau(\eta^i\phi)}^2+\overline{\tau(\eta^i)}^2\tau(\eta^i\phi)^2)^j.
\end{align*}

For $j=0,$ we have 
\begin{align*}
N_j(\eta,n,p)
&=\sum_{1\leqslant i\leqslant \frac{p-1}{2}}\eta^{i}(n^2)=
\frac{p-1}{2}
\end{align*}
if $n\equiv\pm1\bmod p,$ and otherwise
\begin{align*}
N_j(\eta,n,p)
&=\frac{\eta(n^2)(1-\eta(n^2)^{\frac{p-1}{2}})}{1-\eta(n^2)}=\frac{\eta(n^2)(1-\eta(n^{p-1}))}{1-\eta(n^2)}=0
\end{align*}
for $(n,p)=1.$

We henceforth assume $j>0$, for which
we derive that
\begin{align}
N_j(\eta,n,p)
&=\sum_{0\leqslant r\leqslant j}\binom{j}{r}\sum_{1\leqslant i\leqslant \frac{p-1}{2}}\eta^{i}(n^2)(\tau(\eta^i)\overline{\tau(\eta^i\phi)})^{2r}(\overline{\tau(\eta^i)}\tau(\eta^i\phi))^{2(j-r)}\nonumber\\
&=\sum_{0\leqslant r\leqslant j}\binom{j}{r}p^{4(j-r)}\sum_{1\leqslant i\leqslant \frac{p-1}{2}}\eta^{i}(n^2)(\tau(\eta^i)\overline{\tau(\eta^i\phi)})^{4r-2j}+O(p^{2j}).\label{eq:Nj(eta,n,p)}
\end{align}

For $j/2<r\leqslant j,$ Lemma \ref{lm:Gaussmoment-j} yields
\begin{align*}
p^{4(j-r)}\sum_{1\leqslant i\leqslant \frac{p-1}{2}}\eta^{i}(n^2)(\tau(\eta^i)\overline{\tau(\eta^i\phi)})^{4r-2j}\ll (2r-j)p^{2j+\frac{1}{2}}\log p.
\end{align*}
Similarly, for $0\leqslant r<j/2,$ we derive from Lemma \ref{lm:Gaussmoment-j} that
\begin{align*}
&\ \ \ p^{4(j-r)}\sum_{1\leqslant i\leqslant \frac{p-1}{2}}\eta^{i}(n^2)(\tau(\eta^i)\overline{\tau(\eta^i\phi)})^{4r-2j}\\
&=
p^{4r}\sum_{1\leqslant i\leqslant \frac{p-1}{2}}\eta^{i}(n^2)(\overline{\tau(\eta^i)}\tau(\eta^i\phi))^{2j-4r}\\
&\ll (j-2r)p^{2j+\frac{1}{2}}\log p.
\end{align*}
The term with $r=j/2$ will contribute to the main term for $N_j(\eta,n,p)$, which can only occur when $j$ is even.

Combining all the above arguments, we derive for $j>0$ that
\begin{align*}
N_j(\eta,n,p)
&=\frac{\delta_{j}\Delta_n}{2}\binom{j}{j/2}p^{2j+1}+O(j2^jp^{2j+\frac{1}{2}}\log p),
\end{align*}
where $\delta_{j}=1$ or $0$ according to $j$ is even or odd, and $\Delta_n=1$ or $0$ according to $n\equiv\pm1\bmod p$ or not.
Inserting this into \eqref{eq:Mk(p)-decomposition} we find
\begin{align*}
M_\kappa(n,p)
&=2^{\kappa-1}\Delta_np+p^{-2\kappa}\Gamma(\kappa+1)\Delta_n\sum_{j\geqslant 1}\frac{(2p^2)^{\kappa-2j}}{(2j)!\Gamma(\kappa-2j+1)}\binom{2j}{j}p^{4j+1}\\
&\ \ \ \ +O(\kappa4^\kappa p^{\frac{1}{2}}\log p),
\end{align*}
where the first and second terms come from $j=0$ and $j>0$ in \eqref{eq:Mk(p)-decomposition}, respectively.
Then Theorem \ref{thm:K-moment} follows immediately from the observation that
\begin{align*}
1+\sum_{j\geqslant1}\frac{\Gamma(\kappa+1)}{(2j)!\Gamma(\kappa-2j+1)}\binom{2j}{j}2^{-2j}
=\frac{\Gamma(2\kappa+1)}{2^\kappa\Gamma(\kappa+1)^2}
\end{align*}
by Proposition \ref{prop:Gammaidentity} below.

\smallskip

\section{Proof of Theorem \ref{thm:|K|-equidistribution}}\label{sec:|K|-equidistribution-proof}

Theorem \ref{thm:|K|-equidistribution} is an immediate consequence of the method of moments that the convergence in distribution can be detected by the convergence of a sequence of moment sequences. 
See Appendix \ref{sec:methodofmoments} for the precise formulation of this method.

It suffices to prove that, as $p$ tends to infinity along primes, the random variable $\chi\mapsto |K(\chi)|$ converges in distribution to a distribution function with the characteristic function given by
\begin{align*}
f(t)&=\frac{2}{\pi}\int_0^2\frac{\ue^{itx}}{\sqrt{4-x^2}}\ud x.
\end{align*}
In fact, Theorem \ref{thm:K-moment} yields
\begin{align*}
\frac{1}{(p-3)/2}~\sideset{}{^*}\sum_{\substack{\chi\bmod p\\ \chi(-1)=1}}|K(\chi)|^k
&=\frac{\Gamma(k+1)}{\Gamma(k/2+1)^2}+O(2^kp^{-\frac{1}{2}}\log p),
\end{align*}
from which and Lemma \ref{lm:methodofmoments} the characteristic function should be
\begin{align*}
f(t)&=\sum_{k\geqslant0} \frac{\Gamma(k+1)}{\Gamma(k/2+1)^2}\frac{(it)^k}{k!}=\frac{2}{\pi}\int_0^{\frac{\pi}{2}}\sum_{k\geqslant0} \frac{(2it)^k}{k!}(\cos\theta)^k\ud\theta=\frac{2}{\pi}\int_0^{\frac{\pi}{2}}\ue^{2it\cos\theta}\ud\theta
\end{align*}
in view of \eqref{eq:trigonometric-integral}.
This completes the proof of Theorem \ref{thm:|K|-equidistribution} by making a change of variable and Fourier inversion.

\smallskip

\section{Proof of Theorems \ref{thm:K-equidistribution} and \ref{thm:K-mixedmoment}}
\label{sec:K-equidistribution-proof}

We first prove Theorem \ref{thm:K-mixedmoment}.
Without loss of generality we assume $r:=l-k\geqslant 0.$ Recall that $\eta$ is a generator of $\cC_p$, the group of all multiplicative characters mod $p$.
From \eqref{eq:K-Gausssum} we infer
\begin{align*}
M_{k,l}(n,p)
&=\overline{\tau(\phi)}^rp^{-\frac{k+3l}{2}}\sum_{1\leqslant i\leqslant \frac{p-3}{2}}\eta^{i}(4^rn^2)(\tau(\eta^i)\overline{\tau(\eta^i\phi)}+\overline{\tau(\eta^i)}\tau(\eta^i\phi))^{k+l}\\
&=\overline{\tau(\phi)}^rp^{-\frac{k+3l}{2}}\sum_{0\leqslant j\leqslant k+l}\binom{k+l}{j}\sum_{1\leqslant i\leqslant \frac{p-3}{2}}\eta^{i}(4^rn^2)\\
&\ \ \ \times (\tau(\eta^i)\overline{\tau(\eta^i\phi)})^j(\overline{\tau(\eta^i)}\tau(\eta^i\phi))^{k+l-j}.
\end{align*}
If $2j\neq k+l,$ we may derive from Lemma \ref{lm:Gaussmoment-j} that the $i$-sum is at most $O((k+l)p^{k+l+\frac{1}{2}}\log p).$ If $2j=k+l,$ which can only happen when $r$ is even, the $i$-sum is
\begin{align*}
p^{k+l}\sum_{1\leqslant i\leqslant \frac{p-3}{2}}\eta^{i}(4^rn^2)
&=
p^{k+l}\sum_{\chi\in\cC_p^+}\chi(2^rn)+O(p^{k+l})\\
&=
\frac{1}{2}p^{k+l+1}\cdot \Delta_{2^rn}+O(p^{k+l}).
\end{align*}
Combining the above two cases, it follows that
\begin{align*}
M_{k,l}(n,p)
&=\Big(\frac{\overline{\tau(\phi)}}{\sqrt{p}}\Big)^r\binom{k+l}{(k+l)/2}\cdot \frac{\Delta_{2^rn}}{2}\cdot \mathbf{1}_{2\mid r}\cdot p+O((k+l)2^{k+l}p^{\frac{1}{2}}\log p).
\end{align*}
This completes the proof of Theorem \ref{thm:K-mixedmoment} for $l\geqslant k$, and the opposite case can be proven identically.

We now turn to the proof of Theorem \ref{thm:K-equidistribution}. By Lemma \ref{lm:methodofmoments}, it suffices to consider the mixed moment
\begin{align*}
C(k_1,k_2)=\frac{1}{(p-3)/2}~\sideset{}{^*}\sum_{\substack{\chi\bmod p\\ \chi(-1)=1}} (\RRe K(\chi))^{k_1}(\IIm K(\chi))^{k_2}.
\end{align*}
for any non-negative integers $k_1,k_2$.

From the binomial expansion, it follows that
\begin{align*}
C(k_1,k_2)
&=\frac{1}{(p-3)/2}~\sideset{}{^*}\sum_{\chi\bmod p} \Big(\frac{K(\chi)+\overline{K(\chi)}}{2}\Big)^{k_1}\Big(\frac{K(\chi)-\overline{K(\chi)}}{2i}\Big)^{k_2}\\
&=\frac{1}{(p-3)/2}~\sum_{0\leqslant j_1\leqslant k_1}\sum_{0\leqslant j_2\leqslant k_2}\binom{k_1}{j_1}\binom{k_2}{j_2}\frac{(-1)^{k_2-j_2}}{2^{k_1}(2i)^{k_2}}M_{k_1+k_2-j_1-j_2,j_1+j_2}(1,p),
\end{align*}
where $M_{*,*}(1,p)$ is defined by \eqref{eq:mixedmoment}.
Hence
\begin{align*}
C(k_1,k_2)
&=\binom{k_1+k_2}{(k_1+k_2)/2}\mathop{\sum_{0\leqslant j_1\leqslant k_1}\sum_{0\leqslant j_2\leqslant k_2}}_{2^{k_1+k_2}\equiv \pm 2^{2(j_1+j_2)}\bmod p}\binom{k_1}{j_1}\binom{k_2}{j_2}\frac{(-1)^{k_2-j_2}}{2^{k_1}(2i)^{k_2}}\varepsilon_p^{k_1+k_2-2(j_1+j_2)}\\
&\ \ \ +O((k_1+k_2)2^{k_1+k_2}p^{-\frac{1}{2}}\log p),
\end{align*}
where the first term vanishes if $k_1+k_2$ is odd. For any $k_1,k_2=o(\log p)$ (as $p\rightarrow+\infty$), only the terms with $2(j_1+j_2)=k_1+k_2$ survive in the double sum over $j_1,j_2$. Suppose that $k_1+k_2=2k,$ we then have
\begin{align}\label{eq:C(k1,k2)}
C(k_1,k_2)
&=\binom{2k}{k}\mathop{\sum_{0\leqslant j_1\leqslant k_1}\sum_{0\leqslant j_2\leqslant k_2}}_{j_1+j_2=k}\binom{k_1}{j_1}\binom{k_2}{j_2}\frac{(-1)^{k_2-j_2}}{2^{k_1}(2i)^{k_2}}+O(k4^kp^{-\frac{1}{2}}\log p).
\end{align}

Let  $F_p$ be the characteristic function of the joint distribution of $\RRe K(\chi)$ and $\IIm K(\chi)$:
\begin{align*}
F_p(u,v)=\frac{1}{(p-3)/2}~\sideset{}{^*}\sum_{\substack{\chi\bmod p\\ \chi(-1)=1}} \ue^{iu\RRe K(\chi)+iv\IIm K(\chi)}.
\end{align*}
Using the Taypor expansion, we write
\begin{align*}
F_p(u,v)
&=\frac{1}{(p-3)/2}~\sideset{}{^*}\sum_{\substack{\chi\bmod p\\ \chi(-1)=1}} \sum_{k_1\geqslant0}\sum_{k_2\geqslant0}\frac{(iu\RRe K(\chi))^{k_1}(iv\IIm K(\chi))^{k_2}}{k_1!k_2!}\\
&=\sum_{k_1\geqslant0}\sum_{k_2\geqslant0}\frac{(iu)^{k_1}(iv)^{k_2}}{k_1!k_2!}C(k_1,k_2).
\end{align*}
In view of \eqref{eq:C(k1,k2)}, the characteristic function of the limiting distribution, denoted by $F,$ should be
\begin{align}\label{eq:F(u,v)-initial}
F(u,v)
=\mathop{\sum_{k_1\geqslant0}\sum_{k_2\geqslant0}}_{2\mid k_1+k_2}\frac{(iu)^{k_1}(iv)^{k_2}}{k_1!k_2!}\binom{k_1+k_2}{(k_1+k_2)/2}\mathop{\sum_{0\leqslant j_1\leqslant k_1}\sum_{0\leqslant j_2\leqslant k_2}}_{j_1+j_2=(k_1+k_2)/2}\binom{k_1}{j_1}\binom{k_2}{j_2}\frac{(-1)^{k_2-j_2}}{2^{k_1}(2i)^{k_2}},
\end{align}
and it suffices to prove that
\begin{align}\label{eq:F(u,v)-Fourier}
F(u,v)&=\frac{1}{\pi^2}\iint_{x^2+y^2\leqslant4} \frac{\ue^{ir(ux+vy)}\ud x\ud y}{\sqrt{x^2+y^2}\sqrt{4-(x^2+y^2)}}.
\end{align}

Grouping variables $k_1,k_2$ and $j_1,j_2$, and switching summations, we now derive from \eqref{eq:F(u,v)-initial} that
\begin{align*}
F(u,v)
&=\sum_{k\geqslant 0}\binom{2k}{k}\mathop{\sum_{k_1\geqslant0}\sum_{k_2\geqslant0}}_{k_1+k_2=2k}\frac{(iu)^{k_1}(iv)^{k_2}}{k_1!k_2!}\mathop{\sum_{0\leqslant j_1\leqslant k_1}\sum_{0\leqslant j_2\leqslant k_2}}_{j_1+j_2=k}\binom{k_1}{j_1}\binom{k_2}{j_2}\frac{(-1)^{k_2-j_2}}{2^{k_1}(2i)^{k_2}}\\
&=\sum_{k\geqslant 0}\binom{2k}{k}\mathop{\sum_{j_1\geqslant0}\sum_{j_2\geqslant0}}_{j_1+j_2=k}\frac{(-1)^{j_2}}{j_1!j_2!}\mathop{\sum_{k_1\geqslant j_1}\sum_{k_2\geqslant j_2}}_{k_1+k_2=2k}\frac{(\frac{iu}{2})^{k_1}(-\frac{v}{2})^{k_2}}{(k_1-j_1)!(k_2-j_2)!}\\
&=\sum_{k\geqslant 0}\binom{2k}{k}\mathop{\sum_{j_1\geqslant0}\sum_{j_2\geqslant0}}_{j_1+j_2=k}\frac{(\frac{iu}{2})^{j_1}(\frac{v}{2})^{j_2}}{j_1!j_2!}\mathop{\sum_{k_1\geqslant0}\sum_{k_2\geqslant0}}_{k_1+k_2=k}\frac{(\frac{iu}{2})^{k_1}(-\frac{v}{2})^{k_2}}{k_1!k_2!}.
\end{align*}
We then arrive at
\begin{align*}
F(u,v)&=\sum_{k\geqslant 0}\binom{2k}{k}\frac{1}{k!^2}\Big(\frac{iu}{2}+\frac{v}{2}\Big)^k\Big(\frac{iu}{2}-\frac{v}{2}\Big)^k\\
&=\sum_{k\geqslant 0}\binom{2k}{k}\frac{1}{k!^2}\Big(-\frac{u^2+v^2}{4}\Big)^k.
\end{align*}
Appealing to \eqref{eq:trigonometric-integral}, we may express $\binom{2k}{k}$ in terms of integrals, so that
\begin{align*}
F(u,v)
&=\frac{2}{\pi}\int_0^{\frac{\pi}{2}}\sum_{k\geqslant 0}\frac{(-1)^k}{k!^2}(\sqrt{u^2+v^2}\cos\theta)^{2k}\ud\theta\\
&=\frac{2}{\pi}\int_0^{\frac{\pi}{2}}J_0(2\sqrt{u^2+v^2}\cos\theta)\ud\theta,
\end{align*}
where $J_0$ is the Bessel function defined by \eqref{eq:BesselJ-definition} in Appendix \ref{sec:Bessel-Combinatorics}.
Using the integral representation \eqref{eq:BesselJ-integral}, it follows that
\begin{align*}
F(u,v)
&=\frac{1}{\pi^2}\int_0^{\frac{\pi}{2}}\ud\theta\int_{-\pi}^{\pi} e^{2i\sqrt{u^2+v^2}\cos\theta\sin\varphi} \ud\varphi.
\end{align*}
Making changes of variables, we find
\begin{align*}
F(u,v)
&=\frac{1}{\pi^2}\int_0^1\frac{\ud r}{\sqrt{1-r^2}} \int_{-\pi}^{\pi} \ue^{2ir\sqrt{u^2+v^2}\sin \varphi}\ud\varphi\\
&=\frac{1}{\pi^2}\iint_{x^2+y^2\leqslant4} \frac{\ue^{ir(ux+vy)}\ud x\ud y}{\sqrt{x^2+y^2}\sqrt{4-(x^2+y^2)}} \\
\end{align*}
This proves \eqref{eq:F(u,v)-Fourier}, and thus completes the proof of Theorem \ref{thm:K-equidistribution} by Fourier inversion.

\smallskip

\section{Proof of Theorem \ref{thm:K-twistedequidistribution}}\label{sec:K-twistedequidistribution-proof}
The method of moments requires asymptotic evaluations of the moments
\begin{align*}
M_k^*(p):=\sum_{\chi\in\fA(p)}\Big(\frac{\overline{\chi}(2)K(\chi)\tau(\phi)}{\sqrt{p}}\Big)^k
\end{align*}
for any $k\in\bZ^+.$ In fact, we see from \eqref{eq:K-Gausssum} that
\begin{align*}
M_k^*(p)=p^{-k}\sum_{1\leqslant i\leqslant \frac{p-3}{2}}(\tau(\eta^i)\overline{\tau(\eta^i\phi)}+\tau(\eta^i\phi)\overline{\tau(\eta^i)})^k,
\end{align*}
where $\eta$ is a generator of $\cC_p.$
The argument in Section \ref{sec:proof-Theorem-K-moment} yields
\begin{align*}
M_k^*(p)\ll k2^kp^{\frac{1}{2}}\log p
\end{align*}
for each odd $k\geqslant1,$ where the implied constant is absolute.
For each even $k\geqslant1,$ we may appeal to Theorem \ref{thm:K-moment}, getting
\begin{align*}
M_k^*(p)=\frac{1}{2}\binom{k}{k/2}\cdot p+O(k2^kp^{\frac{1}{2}}\log p).
\end{align*}
In view of Lemma \ref{lm:methodofmoments}, the characteristic function $f^*$ of the limiting distribution should be
\begin{align*}
f^*(t)&=\sum_{k\geqslant0} \frac{(2k)!}{k!^2}\frac{(it)^{2k}}{(2k)!}=\sum_{k\geqslant0}\frac{(-t^2)^k}{k!^2}=J_0(2t).
\end{align*}
Thanks to the integral representation \eqref{eq:BesselJ-integral}, we may write
\begin{align*}
f^*(t)&=\frac{1}{\pi}\int_0^{\pi}\ue^{2it\cos\theta}\ud\theta
=\frac{1}{\pi}\int_{-1}^1\frac{\ue^{2itx}}{\sqrt{1-x^2}}\ud x=\frac{1}{\pi}\int_{-2}^2\frac{\ue^{itx}}{\sqrt{4-x^2}}\ud x.
\end{align*}
This proves Theorem \ref{thm:K-twistedequidistribution} by Fourier inversion.

\smallskip

\section{Proof of Corollary \ref{coro:K-moment-L1}}
The proof of Corollary \ref{coro:K-moment-L1} will follow the arguments of Zhang \cite{Zh02}. More precisely, we would like to prove the following preliminary lemma.

\begin{lemma}\label{lm:|L|-expansion}
For each multiplicative character $\chi\bmod p,$ attach a complex coefficient $c(\chi)$ with $|c(\chi)|\leqslant1.$ For all large prime $p,$ we have
\begin{align*}
\sideset{}{^*}\sum_{\chi\bmod p}c(\chi)|L(1,\chi)|
&=\sideset{}{^*}\sum_{\chi\bmod p}c(\chi)|B(\chi,N)|^2+O(N^{-\frac{1}{2}}p+N^{-1}p^{\frac{3}{2}}\log p+p^{\frac{1}{2}}\log p)
\end{align*}
for any $N>1,$ where
\begin{align}\label{eq:B(chi,N)}
B(\chi,N)=\sum_{n\leqslant N}\frac{\chi(n)\beta(n)}{n}.
\end{align}
\end{lemma}

\proof
First, we write
\begin{align*}
B(\chi,N)^2
&=\mathop{\sum\sum}_{n_1,n_2\leqslant N}\frac{\chi(n_1n_2)\beta(n_1)\beta(n_2)}{n_1n_2}=\sum_{n\leqslant N^2}\frac{\chi(n)\gamma(n,N)}{n},
\end{align*}
where
\begin{align*}
\gamma(n,N)&=\mathop{\sum\sum}_{\substack{n_1n_2=n\\ n_1,n_2\leqslant N}}\beta(n_1)\beta(n_2).
\end{align*}
Note that $\gamma(n,N)=(\beta*\beta)(n)=1$ as long as $n\leqslant N.$  This gives
\begin{align*}
B(\chi,N)^2
&=\sum_{n\leqslant N}\frac{\chi(n)}{n}+\sum_{N<n\leqslant N^2}\frac{\chi(n)\gamma(n,N)}{n}.
\end{align*}

On the other hand, for each non-trivial multiplicative character $\chi\bmod p,$ we have
\begin{align*}
L(1,\chi)&=\sum_{n\leqslant N}\frac{\chi(n)}{n}+\int_N^{+\infty}\Big(\sum_{N<n\leqslant y}\chi(n)\Big)\frac{\ud y}{y^2}=\sum_{n\leqslant N}\frac{\chi(n)}{n}+O(N^{-1}p^{\frac{1}{2}}\log p)
\end{align*}
using the P\'olya--Vinogradov inequality for incomplete character sums.

In this way, we derive that
\begin{align*}
&\ \ \ \ \ \sideset{}{^*}\sum_{\chi\bmod p}c(\chi)|L(1,\chi)|-\sideset{}{^*}\sum_{\chi\bmod p}c(\chi)|B(\chi,N)|^2\\
&\ll\sideset{}{^*}\sum_{\chi\bmod p}\Big|\sum_{N<n\leqslant N^2}\frac{\chi(n)\gamma(n,N)}{n}\Big|+N^{-1}p^{\frac{3}{2}}\log p=\Sigma+N^{-1}p^{\frac{3}{2}}\log p,
\end{align*}
say. By Cauchy--Schwarz,
\begin{align*}
\Sigma^2
&\ll p\sideset{}{^*}\sum_{\chi\bmod p}\Big|\sum_{N<n\leqslant N^2}\frac{\chi(n)\gamma(n,N)}{n}\Big|^2\\
&\ll p^2\mathop{\sum\sum}_{\substack{N<n_1,n_2\leqslant N^2\\ n_1\equiv n_2\bmod p}}\frac{1}{n_1n_2}\\
&\ll p^2\sum_{N<n\leqslant N^2}\sum_{0\leqslant |l|\leqslant N^2/p}\frac{1}{n(n+lp)}\\
&\ll N^{-1}p^2+p\log^2p,
\end{align*}
where we have used the trivial bound $\|\beta\|_\infty\leqslant1.$
The lemma then follows by combining all above estimates.
\endproof

We are now ready to prove Corollary \ref{coro:K-moment-L1}. Taking $c(\chi)=|K(\chi)/2|^{2\kappa}$ and $N>p$ in Lemma \ref{lm:|L|-expansion}, we have
\begin{align*}
\sideset{}{^*}\sum_{\chi\bmod p}|K(\chi)|^{2\kappa}|L(1,\chi)|
&=\sideset{}{^*}\sum_{\chi\bmod p}|K(\chi)|^{2\kappa}|B(\chi,N)|^2+O(p^{\frac{1}{2}}\log p).
\end{align*} 
Invoking the expression \eqref{eq:B(chi,N)}, it follows that
\begin{align*}
\sideset{}{^*}\sum_{\chi\bmod p}|K(\chi)|^{2\kappa}|L(1,\chi)|
&=\mathop{\sum\sum}_{\substack{m,n\leqslant N\\(mn,p)=1}}\frac{\beta(m)\beta(n)}{mn}\sideset{}{^*}\sum_{\chi\bmod p}\chi(m\overline{n})|K(\chi)|^{2\kappa}+O(p^{\frac{1}{2}}\log p).
\end{align*} 
Then Theorem \ref{thm:K-moment} yields
\begin{align*}
\sideset{}{^*}\sum_{\chi\bmod p}|K(\chi)|^{2\kappa}|L(1,\chi)|
&=\frac{\Gamma(2\kappa+1)}{2\Gamma(\kappa+1)^2}\mathop{\sum\sum}_{\substack{m,n\leqslant N\\m\equiv\pm n\bmod p\\ (mn,p)=1}}\frac{\beta(m)\beta(n)}{mn}\cdot p+O(4^{\kappa}p^{\frac{1}{2}}\log p).
\end{align*} 
In fact,
\begin{align*}
\mathop{\sum\sum}_{\substack{m,n\leqslant N\\m\equiv n\bmod p\\ (mn,p)=1}}\frac{\beta(m)\beta(n)}{mn}
&=\sum_{\substack{n\leqslant N\\(n,p)=1}}\frac{\beta(n)^2}{n^2}+O(p^{-1})
=\fS+O(p^{-1}),
\end{align*} 
and 
\begin{align*}
\mathop{\sum\sum}_{\substack{m,n\leqslant N\\m\equiv -n\bmod p\\ (mn,p)=1}}\frac{\beta(m)\beta(n)}{mn}
&=O(p^{-1}).
\end{align*} 
Now Corollary \ref{coro:K-moment-L1} follows immediately.

\smallskip

\section{Proof of Theorem \ref{thm:K-moment-L1/2}}
We start the proof of Theorem \ref{thm:K-moment-L1/2} from the following approximate functional equation; see \cite[Theorem 5.3]{IK04}.
\begin{lemma}\label{lm:AFE}
Let $\chi$ be a multiplicative character mod $p.$
Let $G$ be an entire function satisfying the decay condition $G(s)\ll(1+|s|)^{-A}$ for any $A > 0$ in any strip $|\Re(s)|<B$. Furthermore assume $G(s)=G(-s)$ and that $G$ has a double zero at each $s\in\frac{1}{2}+\bZ.$ Then we have
\begin{align*}
|L(1/2,\chi)|^2&=2\mathop{\sum\sum}_{m,n\geqslant1}\frac{\chi(m\overline{n})}{\sqrt{mn}}V_\fa\Big(\frac{mn}{p}\Big),
\end{align*}
where $\fa=(1-\chi(-1))/2$ and
\begin{align*}
V_\fa(x)=\frac{1}{2\pi i}\int_{(1)}\pi^{-s}\frac{G(s)}{s}\frac{\Gamma(\frac{\frac{1}{2}+s+\fa}{2})^2}{\Gamma(\frac{\frac{1}{2}+\fa}{2})^2}x^{-s}\ud s.
\end{align*}
In particular, for all $\ell\geqslant0$ we have
\begin{align*}
x^\ell V_\fa^{(\ell)}(x)\ll (1+|x|)^{-A}
\end{align*}
for any $A\geqslant0$
and 
\begin{align*}
x^\ell V_\fa^{(\ell)}(x)=\delta_\ell+O(|x|^{\frac{1}{2}-\varepsilon}),
\end{align*}
where $\delta_\ell=1$ or $0$ according to $\ell=0$ or $\ell>0,$ and the implied constants depends on $A$ and $\varepsilon.$
\end{lemma}
\begin{remark}
A typical example of choosing $G$ would be $G(s)=e^{s^2}\cos^2(\pi s).$
\end{remark}

By virtue of Lemma \ref{lm:AFE} we write
\begin{align*}
\sideset{}{^*}\sum_{\chi\bmod p}|K(\chi)|^{2\kappa}|L(1/2,\chi)|^2&=2\mathop{\sum\sum}_{m,n\geqslant1}\frac{1}{\sqrt{mn}}V_0\Big(\frac{mn}{p}\Big)M_\kappa(m\overline{n},p).
\end{align*}
The main term comes from $m\equiv\pm n\bmod p$. It then follows that
\begin{align}\label{eq:momentL1/2-decomposition}
\sideset{}{^*}\sum_{\chi\bmod p}|K(\chi)|^{2\kappa}|L(1/2,\chi)|^2&=2M_\kappa(1,p)\cM_\kappa(p)+\cE_\kappa(p),
\end{align}
where
\begin{align*}
\cM_\kappa(p)=\mathop{\sum\sum}_{\substack{m,n\geqslant1\\ m\equiv\pm n\bmod p\\(mn,p)=1}}\frac{1}{\sqrt{mn}}V_0\Big(\frac{mn}{p}\Big),
\end{align*}
\begin{align*}
\cE_\kappa(p)&=2\mathop{\sum\sum}_{\substack{m,n\geqslant1\\ m\not\equiv\pm n\bmod p\\(mn,p)=1}}\frac{1}{\sqrt{mn}}V_0\Big(\frac{mn}{p}\Big)M_\kappa(m\overline{n},p).
\end{align*}

In fact,
\begin{align}\label{eq:Mkappa(p):mainterm}
\cM_\kappa(p)
&=\frac{1}{p-1}\sideset{}{^*}\sum_{\substack{\chi\bmod p\\ \chi(-1)=1}}|L(1/2,\chi)|^2.
\end{align}
Note that
\begin{align*}
\mathop{\sum\sum}_{\substack{m,n\geqslant1\\ mn>p^{1+\varepsilon}\\ m\not\equiv\pm n\bmod p\\(mn,p)=1}}\frac{1}{\sqrt{mn}}V_0\Big(\frac{mn}{p}\Big)M_\kappa(m\overline{n},p)
&\ll 4^\kappa p^{\frac{1}{2}}\log p\sum_{n>p^{1+\varepsilon}}\frac{\tau(n)}{\sqrt{n}}\Big|V_0\Big(\frac{n}{p}\Big)\Big|\ll 4^\kappa
\end{align*}
due to the rapid decay of $V_0.$ Hence
\begin{align*}
\cE_\kappa(p)&=2\mathop{\sum\sum}_{\substack{m,n\geqslant1\\ mn\leqslant p^{1+\varepsilon}\\ m\not\equiv\pm n\bmod p\\(mn,p)=1}}\frac{1}{\sqrt{mn}}V_0\Big(\frac{mn}{p}\Big)M_\kappa(m\overline{n},p)+O(4^\kappa).
\end{align*}
Following the similar arguments in proving Lemma \ref{lm:Gaussmoment-j}, we need to study
\begin{align*}
\fS_k(h,p):=\mathop{\sum\sum}_{\substack{m,n\geqslant1\\ mn\leqslant p^{1+\varepsilon}\\ (mn,p)=1}}\frac{1}{\sqrt{mn}}V_0\Big(\frac{mn}{p}\Big)(1-\phi(g^hmn))\kl_k(g^h(\overline{m}n)^2,p)
\end{align*}
for each fixed integer $k\geqslant2$. Trivially, we have
\begin{align*}
\fS_k(h,p)\ll kp^{\frac{1}{2}+\varepsilon},
\end{align*}
and the goal here is to prove
\begin{align}\label{eq:Sk(h,p)-upperbound}
\fS_k(h,p)\ll k^Bp^{\frac{1}{2}-\delta}
\end{align}
for some positive $\delta>0$ and $B>0$ uniformly in $h$, from which we then conclude that
\begin{align*}
\cE_\kappa(p)&\ll 5^\kappa p^{1-\delta+\varepsilon}.
\end{align*}
Combining this with \eqref{eq:momentL1/2-decomposition} and \eqref{eq:Mk(p)-decomposition}, we may conclude the following proposition.

\begin{proposition}\label{prop:K-moment-L1/2}
Let $\kappa$ be a fixed positive real number. For all large primes $p,$ we have
\begin{align*}
\sideset{}{^*}\sum_{\chi\bmod p}|K(\chi)|^{2\kappa}|L(1/2,\chi)|^2=\frac{2M_\kappa(1,p)}{p-1}\sideset{}{^*}\sum_{\substack{\chi\bmod p\\ \chi(-1)=1}}|L(1/2,\chi)|^2+O(5^\kappa p^{1-\delta+\varepsilon}),
\end{align*}
where $\delta$ is determined by $\eqref{eq:Sk(h,p)-upperbound}$ and the implied constant depends only on $\varepsilon.$
\end{proposition}

The second moment of $L(1/2,\chi)$ with $\chi$ varying over even primitive characters mod $p$ is already known to many authors with suitable generalizations or applications, see \cite{IS99}, \cite{Co07} and \cite{Yo11} for instance. The following version is sufficient for our purpose.

\begin{lemma}\label{lm:L1/2-secondmoent}
We have
\begin{align*}
\sideset{}{^*}\sum_{\chi\bmod p}|L(1/2,\chi)|^2=\frac{(p-1)^2}{2p}(\log\frac{p}{8\pi}-\frac{\pi}{2}+\gamma)+O(\sqrt{p}),
\end{align*}
where the implied constant is absolute.
\end{lemma}

Now Theorem \ref{thm:K-moment-L1/2} follows by combining Lemma \ref{lm:L1/2-secondmoent} and Proposition \ref{prop:K-moment-L1/2} with some admissible $\delta>0.$

It remains to prove the claim in \eqref{eq:Sk(h,p)-upperbound} for any $\delta<1/8$.
The key point will be controlling the cancellations while summing over $m$ and $n$. To this end, we may invoke a general result of Fouvry, Kowalski and Michel on (smooth) bilinear forms for general trace functions, see \cite[Theorem 1.16]{FKM14}. 

To state the work of Fouvry, Kowalski and Michel, we assume $K$ is an isotypic trace weight associated to the $\ell$-adic sheaf $\cF$ modulo $p$ of conductor $\fc(K)$. See \cite[Definition 1.3]{FKM14} for the precise definitions.

\begin{lemma}\label{lm:FKM}
 Let $M,N, X\geqslant 1$ be parameters with
$X/4\leqslant MN\leqslant X$. Let $W_1,W_2,W_3$ be smooth functions with supports in $[1,2]$ satisfying the growth condition
\begin{align}\label{eq:W-decay}
x^jW_\ell^{(j)}(x)\ll_j Q_\ell^j
\end{align}
for some $Q_\ell\geqslant1$ with $\ell=1,2,3$ and all $j\geqslant0.$
Then we have
\begin{align*}
\mathop{\sum_m\sum_n}K(mn)\Big(\frac{m}n\Big)^{it}
& W_1\Big(\frac{m}M\Bigr)W_2\Bigl(\frac{n}N\Bigr) W_3\Bigl(\frac{mn}{X}\Big)\\
&\ll (\fc(K)+|t|+Q_1+Q_2)^B Q_3
X\Big(1+\frac{p}{X}\Big)^{1/2}p^{-\eta}
\end{align*}
for $t\in\bR$ and for any $\eta< 1/8$ and some constant $B\geqslant 1$
depending on $\eta$ only. The implicit constant depends only on $\eta$ and the implied constants in $\eqref{eq:W-decay}.$
\end{lemma}

Lemma \ref{lm:FKM} can not be applied directly to the bilinear form $\fS_k(h,p)$ since one has the trace function $(m,n)\mapsto \kl_k(g^h(\overline{m}n)^2).$ In fact, we would like to apply Poisson summation to one variable initially. To begin with, we write
\begin{align*}
\fS_k(h,p)=\fS_k'(h,p)-\phi(g^h)\fS_k''(h,p)
\end{align*}
with
\begin{align*}
\fS_k'(h,p):=\mathop{\sum\sum}_{\substack{m,n\geqslant1\\ mn\leqslant p^{1+\varepsilon}\\ (mn,p)=1}}\frac{1}{\sqrt{mn}}V_0\Big(\frac{mn}{p}\Big)\kl_k(g^h(\overline{m}n)^2,p),
\end{align*}
\begin{align*}
\fS_k''(h,p):=\mathop{\sum\sum}_{\substack{m,n\geqslant1\\ mn\leqslant p^{1+\varepsilon}\\ (mn,p)=1}}\frac{1}{\sqrt{mn}}V_0\Big(\frac{mn}{p}\Big)\phi(mn)\kl_k(g^h(\overline{m}n)^2,p).
\end{align*}
Due to the similar features between $\fS_k'(h,p)$ and $\fS_k''(h,p)$, we only give the details for $\fS_k''(h,p)$, and in order to ease the presentation, we consider the following smoothed version
\begin{align*}
\fS_k:=\mathop{\sum_m\sum_n}U\Big(\frac{m}{M}\Big)V\Big(\frac{n}{N}\Big)\phi(mn)\kl_k(g^h(\overline{m}n)^2,p),
\end{align*}
where $N\geqslant M,$ $MN\leqslant p^{1+\varepsilon}$ and $U,V$ are smooth functions with compact supports in $[1,2]$ satisfying 
\[x^jU^{(j)}(x)\ll_j 1,\ \ \ x^jV^{(j)}(x)\ll_j 1.\]
We would like to show that
\begin{align}\label{eq:Sk-upperbound}
\fS_k\ll k^Bp^{1-\delta},
\end{align}
which would yield \eqref{eq:Sk(h,p)-upperbound} with additional analysis (partition of unity, for instance).

Applying Poisson to the $n$-sum in $\fS_k,$ we obtain
\begin{align*}
\fS_k&=\frac{N}{\sqrt{p}}\mathop{\sum_m\sum_n}U\Big(\frac{m}{M}\Big)\widehat{V}\Big(\frac{nN}{p}\Big)\phi(m)\cL_k(m,n;p,g^h),
\end{align*}
where
\begin{align*}
\cL_k(m,n;p,a):=\frac{1}{\sqrt{p}}\sum_{x\bmod p}\phi(x)\kl_k(a(\overline{m}x)^2,p)\ue\Big(\frac{-nx}{p}\Big).
\end{align*}
Making the change of variable $x\mapsto mx$, we have
\begin{align*}
\cL_k(m,n;p,a)&=\frac{\phi(m)}{\sqrt{p}}\sum_{x\bmod p}\phi(x)\kl_k(ax^2,p)\ue\Big(\frac{-mnx}{p}\Big)
\end{align*}
In this way, we may write
\begin{align*}
\fS_k&=\frac{N}{\sqrt{p}}\mathop{\sum_{(m,p)=1}\sum_n}U\Big(\frac{m}{M}\Big)\widehat{V}\Big(\frac{nN}{p}\Big)F_k(mn;p,g^h),
\end{align*}
where
\begin{align*}
F_k(y;p,a):=\frac{1}{\sqrt{p}}\sum_{x\bmod p}\phi(x)\kl_k(ax^2,p)\ue\Big(\frac{-yx}{p}\Big).
\end{align*}

As a function in $y,$ $F_k(y;p,a)$ is the (normalized) Fourier transform of $\phi(x)\kl_k(ax^2,p).$ According to Laumon \cite{La87} on $\ell$-adic Fourier transforms, the underlying sheaf related to $y\mapsto F_k(y;p,a)$, for each fixed $(a,p)=1,$ is 
 geometrically irreducible and isotypic since $\cL_\phi\otimes[x\mapsto ax^2]^*\cK\ell_k$ is exactly so, where $\cL_\phi$ denotes the Kummer sheaf given by $\phi$ and $\cK\ell_k$ is the hyper-Kloosterman sheaf of rank $k$. Moreover, the computation of Fouvry, Kowalski and Michel \cite{FKM15} shows that the conductor of the trace function $y\mapsto F_k(y;p,a)$ is at most $10k^2.$
 
 We are now in a good position to apply Lemma \ref{lm:FKM}, getting 
 \begin{align*}
\fS_k&\ll k^B p^{1-\eta}
\end{align*}
for any $\eta<1/8.$ This proves  \eqref{eq:Sk-upperbound} and thus  \eqref{eq:Sk(h,p)-upperbound} with any $\delta<1/8$.

\smallskip

\section{Concluding remarks}\label{sec:remarks}

\subsection{Gauss sums with quadratic phases}
As an extension to $\tau(n,\chi)$, we define Gauss sums with general monomial phases in additive characters:
\begin{align*}
\tau_k(n,\chi)=\sum_{v\bmod p}\chi(v)\ue\Big(\frac{nv^k}{p}\Big).
\end{align*}
In particular, we have $\tau_1(n,\chi)=\tau(n,\chi).$ The case $k=2$ is ultimately related to $K(\chi)$ as mentioned in the first section. In fact, if $\chi$ is non-trivial, we have
\begin{align*}
|\tau_2(n,\chi)|^2
&=\mathop{\sideset{}{^*}\sum\sideset{}{^*}\sum}_{u,v\bmod p}\chi(u)\overline{\chi}(v)\ue\Big(\frac{n(u^2-v^2)}{p}\Big)\\
&=\mathop{\sideset{}{^*}\sum\sideset{}{^*}\sum}_{u,v\bmod p}\chi(u)\ue\Big(\frac{nv^2(u^2-1)}{p}\Big)\\
&=\sideset{}{^*}\sum_{u\bmod p}\chi(u)\sum_{v\bmod p}\ue\Big(\frac{nv^2(u^2-1)}{p}\Big)\\
&=(1+\chi(-1))p+\sideset{}{^*}\sum_{\substack{u\bmod p\\ u\not\equiv\pm1\bmod p}}\chi(u)\tau(n(u^2-1),\phi)\\
&=(1+\chi(-1))p+\phi(n)\tau(\phi)\sideset{}{^*}\sum_{u\bmod p}\chi(u)\phi(u^2-1).
\end{align*}
From Lemma \ref{lm:K-transform} it follows that
\begin{align*}
|\tau_2(n,\chi)|^2
&=(1+\chi(-1)+\varepsilon_p\phi(n)\overline{\chi}(2)K(\chi))p.
\end{align*}
and thus all the above moments related to $K(\chi)$ can be reformulated in terms of $|\tau_2(n,\chi)|.$
We omit the details here, and the statements with proofs are left to curious readers.

\subsection{High-dimensional generalization of $K(\chi)$}
As mentioned in the first section, it is very possible to extend our arguments to the general sum $K_k(a,\chi)$ as defined by \eqref{eq:K-hyper}. We will show that $K_k(a,\chi)$ can be interpreted in terms of Gauss sums and the following high-dimensional Jacobi sums:
\begin{align}\label{eq:Jacobisum-hyper}
J(\chi_1,\chi_2,\cdots,\chi_k)=\mathop{\sum\cdots\sum}_{\substack{v_1,v_2,\cdots,v_k\bmod p\\ v_1+v_2+\cdots+v_k\equiv 1\bmod p}}\chi_1(v_1)\chi_2(v_2)\cdots \chi_k(v_k),
\end{align}
where $\chi_1,\chi_2,\cdots,\chi_k$ can be any multiplicative characters mod $p$.

Note that
\begin{align*}
K_k(a,\chi)
&=p^{\frac{1-k}{2}}\sideset{}{^*}\sum_{x\bmod p}\chi(x)N_k(a,x),
\end{align*}
where $N_k(a,x)$ counts the number of tuples $(a_1,a_2,\cdots,a_k)\in[1,p]^k$ satisfying
\begin{align*}
a_1a_2\cdots a_k\equiv a\bmod p,\ \ a_1+a_2+\cdots+a_k\equiv x\bmod p.
\end{align*}
Using orthogonality of multiplicative characters, we may write
\begin{align*}
N_k(a,x)&=\frac{1}{\varphi(p)}\sum_{\rho\bmod p}\overline{\rho}(a)\mathop{\sum\cdots\sum}_{\substack{a_1,a_2,\cdots,a_k\bmod p\\ a_1+a_2+\cdots+a_k\equiv x\bmod p}}\rho(a_1)\rho(a_2)\cdots\rho(a_k)\\
&=\frac{1}{\varphi(p)}\sum_{\rho\bmod p}\overline{\rho}(a)\rho(x)^kJ(\rho,\rho,\cdots,\rho)
\end{align*}
for all $(x,p)=1,$ where there are $k$ copies of $\rho$ in $J(\cdots).$ It then follows that
\begin{align*}
K_k(a,\chi)
&=p^{\frac{1-k}{2}}\sum_{\substack{\rho\bmod p\\ \rho^k=\chi}}\rho(a)\overline{J(\rho,\rho,\cdots,\rho)}.
\end{align*}
We may generalize Lemma \ref{lm:Jacobi-Gauss} to the situation of high-dimensional Jacobi sums, so that
\begin{align*}
K_k(a,\chi)
&=p^{\frac{-1-k}{2}}\tau(\chi)\sum_{\substack{\rho\bmod p\\ \eta^k=\chi}}\rho(a)\overline{\tau(\rho)}^k
\end{align*}
if $\chi$ is non-trivial, where $\tau(\rho)$ is a Gauss sum defined as before. This expression allows one to study the distribution and moments of $K_k(a,\chi)$ by virtue of Gauss sums.

\smallskip

\appendix

\section{Special functions and a combinatorial identity}
\label{sec:Bessel-Combinatorics}
\subsection{Special functions}
Bessel functions were originally introduced to solve a differential equation of second order. We now formulate the Bessel function of the first kind
\begin{align}\label{eq:BesselJ-definition}
J_\nu(z)=\sum_{k\geqslant0}\frac{(-1)^k(z/2)^{2 k+\nu}}{k!\Gamma(\nu+k+1)}
\end{align}
for $|\arg z|<\pi$ and  $\nu\in\bC.$ 
We also need integral representations for $J_\nu.$ In particular, for each non-negative integer $n$, the definition \eqref{eq:BesselJ-definition} is valid for all $z\in\bC,$ for which we have
\begin{align}\label{eq:BesselJ-integral}
J_{n}(z) &=\frac{1}{2 \pi} \int_{-\pi}^{\pi} \ue^{-ni\varphi+iz\sin\varphi}\ud\varphi.
\end{align}
See \cite[Formula 1 in Section 8.41, pp.912]{GR07}.

We also need the following identity
\begin{align}\label{eq:trigonometric-integral}
\int_0^{\frac{\pi}{2}} (\cos\theta)^{2\mu}\ud\theta=\frac{\pi}{2^{2\mu+1}}\frac{\Gamma(2\mu+1)}{\Gamma(\mu+1)^2},\ \ \RRe\mu>-\tfrac{1}{2}
\end{align}
which is used in two directions in this paper; see \cite[Formula 9 in Section 3.63, pp.397]{GR07}.

\subsection{A combinatorial identity}

We now present a combinatorial identity involving Gamma functions, which might have appeared in some existing references. 

\begin{proposition}\label{prop:Gammaidentity}
Let $\kappa\geqslant0$ be a fixed real number. Then we have
\begin{align*}
\sum_{j\geqslant0}\frac{\Gamma(\kappa+1)}{(2j)!\Gamma(\kappa-2j+1)}\binom{2j}{j}2^{-2j}=\frac{\Gamma(2\kappa+1)}{2^\kappa\Gamma(\kappa+1)^2},
\end{align*}
\end{proposition}

The proof of Proposition \ref{prop:Gammaidentity} is a routine application of the analytic theory of Gamma and Bessel functions. Denote by $I$ the series in question. We start from the following contour representation:
\begin{align*}
\frac{1}{\Gamma(z)}=\frac{1}{2\pi i}\int_\cH s^{-z}\ue^s\ud s,
\end{align*}
where $\cH=\cH(r)$ is the Hankel contour shown as follows (suppose the radius of the circle around $0$ is $r$):
\begin{center}
\includegraphics[scale=0.23]{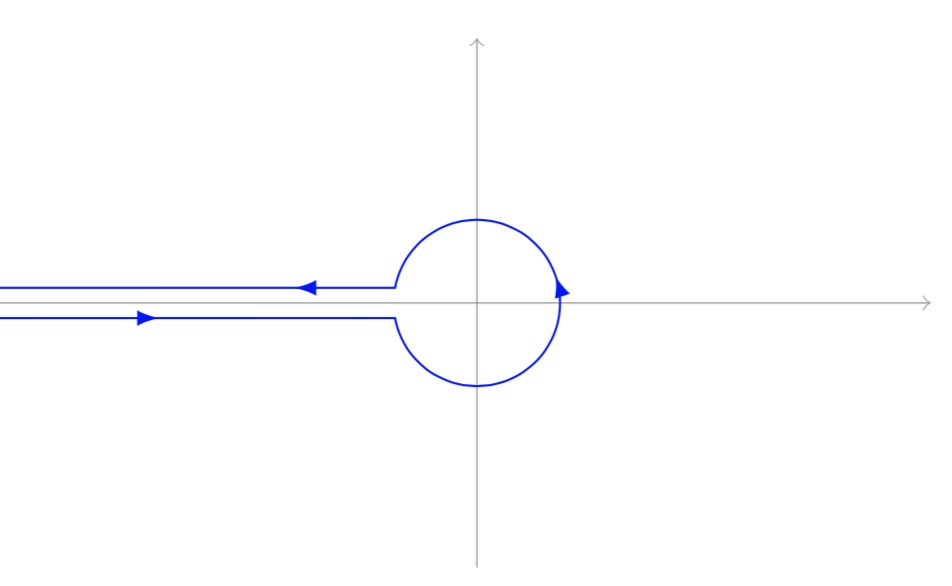}
\end{center}
We now have
\begin{align*}
I&=\sum_{j\geqslant0}\frac{\Gamma(\kappa+1)}{(2j)!}\binom{2j}{j}2^{-2j}\frac{1}{2\pi i}\int_\cH s^{-(\kappa-2j+1)}\ue^s\ud s\\
&=\Gamma(\kappa+1)\frac{1}{2\pi i}\int_\cH s^{-(\kappa+1)}\ue^s\ud s\sum_{j\geqslant0}\frac{(s^2/4)^j}{(j!)^2}\\
&=\Gamma(\kappa+1)\frac{1}{2\pi i}\int_\cH s^{-(\kappa+1)}\ue^sJ_0(is)\ud s,
\end{align*}
where $J_\nu(s)$ is the Bessel function defined by \eqref{eq:BesselJ-definition}. From the integral representation \eqref{eq:BesselJ-integral}, it follows that
\begin{align}\label{eq:I-I(kappa,theta)}
I&=\frac{\Gamma(\kappa+1)}{2\pi}\int_{-\pi}^\pi I(\kappa,\theta)\ud\theta,
\end{align}
where
\begin{align*}
I(\kappa,\theta):=\frac{1}{2\pi i}\int_\cH s^{-(\kappa+1)}\ue^{s(1-\sin\theta)}\ud s.
\end{align*}
For any $\theta\not\equiv\pi/2\bmod\pi,$ we always have $1-\sin\theta>0.$ Making the change of variable $s\mapsto s/(1-\sin\theta),$ we infer
\begin{align*}
I(\kappa,\theta)&=\frac{(1-\sin\theta)^\kappa}{2\pi i}\int_{\cH(r(1-\sin\theta))} s^{-(\kappa+1)}\ue^{s}\ud s.
\end{align*}
Again by the Hankel contour representation, we find
\begin{align*}
I(\kappa,\theta)&=\frac{(1-\sin\theta)^\kappa}{\Gamma(\kappa+1)},
\end{align*}
from which and \eqref{eq:I-I(kappa,theta)} it follows that
\begin{align*}
I&=\frac{1}{2\pi}\int_{-\pi}^\pi (1-\sin\theta)^\kappa\ud\theta=\frac{1}{2\pi}\int_{-\pi}^\pi (1+\cos\theta)^\kappa\ud\theta=\frac{2^{\kappa+1}}{\pi}\int_0^{\frac{\pi}{2}} (\cos\theta)^{2\kappa}\ud\theta.
\end{align*}
Now Proposition \ref{prop:Gammaidentity} follows immediately from the identity \eqref{eq:trigonometric-integral}.

\smallskip

\section{The method of moments}\label{sec:methodofmoments}
In probability theory, the method of moments can be utilized to prove convergence in distribution by showing the convergence of a sequence of moment sequences.
The following result was proven independently by Wintner \cite{Wi28} and Fr\'echet and Shohat \cite{FS31}.
Its prototype is due to Markov, who would like to call this the {\it Second Limit Theorem}. The following formulation can also found in many references; see Billingsley \cite[Theorem 30.1, Exercise 30.6]{Bi95} for instance.

\begin{lemma}\label{lm:methodofmoments}
Let $n$ be a positive integer.
Let $\mu$ be a probability measure on $\bR^n$ such that
\begin{align}\label{eq:methodofmoments-assumption1}
\int_{\bR^n} |x_j|^k\mu(\ud x)<+\infty
\end{align}
for all $1\leqslant j\leqslant n$ and $k\geqslant1.$ Suppose the power series in $\theta$
\begin{align}\label{eq:methodofmoments-assumption2}
\sum_{k\geqslant0}\frac{\theta^k}{k!}\int_{\bR^n} |x_j|^k\mu(\ud x)
\end{align}
has a positive radius $r$ of convergence for each $1\leqslant j\leqslant n$. Then $\mu$ is the only probability measure with the moments 
\begin{align*}
m(k_1,\cdots,k_n)=\int_{\bR^n} x_1^{k_1}\cdots x_n^{k_n} \mu(\ud x)
\end{align*}
for non-negative integers $k_1,\cdots,k_n,$ and the characteristic function $f$ of this distribution has the representation 
\begin{align*}
f(t_1,\cdots,t_n)=\mathop{\sum\cdots\sum}_{k_1,\cdots,k_n\geqslant0}m(k_1,\cdots,k_n)\frac{(it_1)^{k_1}\cdots (it_n)^{k_n}}{k_1!\cdots k_n!}
\end{align*}
 in the disc $\cD^n$ with $\cD=\{t\in\bC:|t|<r\}.$
 
Moreover, if there is a sequence of probability measures $\{\mu_\ell\}_{\ell\geqslant1}$ satisfying the corresponding assumptions $\eqref{eq:methodofmoments-assumption1}$ and $\eqref{eq:methodofmoments-assumption2}$, and suppose that
the moments satisfy
\begin{align*}
\lim_{\ell\rightarrow+\infty}\int_{\bR^n} x_1^{k_1}\cdots x_n^{k_n} \mu_\ell(\ud x)=m(k_1,\cdots,k_n),
\end{align*}
then $\mu_\ell$ converges weakly to $\mu.$
\end{lemma}

\bibliographystyle{plainnat}

\bigskip

\end{document}